# Photovoltaic Generation in Distribution Networks: Optimal vs. Random Installation

Hamidreza Sadeghian, Zhifang Wang
Department of Electrical and Computer Engineering
Virginia Commonwealth University
Richmond, VA, USA
E-mail: {sadeghianh, zfwang}@vcu.edu

*Abstract*— Nowadays common practice in deploying photovoltaic distributed generations (PVDGs) is customer-based installation in the distribution network. Increasing level of PVDG applications and expedite approval by utilities have raised concern about the negative impacts of PVDG installations on the distribution network operations such as reverse power flows and undesirable voltage fluctuations. One potential solutions is to optimize the siting and sizing of these distributed renewable generation resources. This paper presents a comparative study on both optimal and randomized installation of PVDGs with the latter modeling real life customer-based renewable integration. The proposed models examine and compare the impacts of PVDG installation on distribution network operation. Numerical simulations have been performed on a local distribution network model with realistic load profiles, GIS information, local solar insolation, and feeder and voltage settings. It is found that when the distribution system has a medium penetration ratio optimal PVDG installations may introduce essential improvements in terms of voltage deviation and energy loss reduction than randomized installation. However, if the penetration ratio is very low or extremely high there will be not significant difference between the two.

## I. INTRODUCTION

In USA, although in few cases the DG application process requires evaluation of all interconnection requests, all utilities offer expedite approval for small scale DGs (e.g., 25 kW and below) requested by customers[1]. The increasing number of applications raised concerns of many distribution network operators (DNOs) that they do not feel able to guarantee reliability and quality to other customers once they allow large aggregation of DGs to be connected to the distribution network [1]. From an electrical perspective, PVDGs are sources of electrical energy at distribution networks for which these networks have not been designed initially and allocating PV distributed generation (PVDG) systems in distribution systems may inflict unwanted challenges in traditional power systems, which have been designed radial and unidirectional [2]. The most common potential concerns caused by solar power are steady-state overvoltage, impacts on system losses, and issues with voltage regulating devices, protection, and voltage fluctuation [3].

Optimal sizing and sitting of DGs as a solution to address the DG impacts on the electrical network have been extensively studied in the past few decades [4]. Examples of studies include the assessment of maximum DG penetration ratios [5], [6], rooftop PVDG on residential customers [7], [8], numerous analytical [9], [10] probabilistic [11], [12], and heuristic approaches [13], [14] for DG sitting and sizing, aiming to fulfill different technical and/or economical criteria. In [15] authors compared the centralized utility-based DGs in which the utility owns and operates the renewable DGs with the decentralized consumer-based DGs in which each consumer owns and operates the renewable DGs. For detailed reviews see Refs. [16], [17]. The review of the literature shows that to date, numerous studies related to PVDG in distribution networks have focused on optimal installation of PVDGs and mitigating high PV penetration issues. However, there is a gap between proposed studies by researchers and common practice for PVDG installations in real world. Considering customer decision based PVDG applications, it is important to model customer behavior to study different PVDG penetration impacts and, most importantly, compare with optimal allocation to have a comparative insight on both situations. In other words, this study is intended to bridge randomized and optimized PVDG installation and elucidate what will happen if customers are allowed to freely install rooftop PVDGs on their premises and is it worthy to optimally install PVDG systems?

The main contributions of this paper can be listed as: a) Introducing a deployment framework that allows optimization in both the size and location of PVDG to minimize energy loss and voltage deviation subjected to distributed PV constraints and operational constraints; b) Development of a model to mimic customer behavior in PVDG deployment and investigated the impacts of randomized PVDG installation; and c) Completing a comparative study on both optimal installation and free customer decision based PVDG installation. The combination of an optimized model with a statistical model to study a PVDG installation impacts, allows to identify the outcomes for different PVDG deployment policies towards possible strategies to maximize advantages and minimize negative impacts of PVDGs.

The rest of the paper is organized as follows: Section II describes the modeling approach for a distribution network. Section III explains the problem formulation and methodology of comparative analyses. Simulation results are presented in section IV. Finally, conclusions are drawn in section V.

## II. Electrical distribution network modeling

A radial electrical distribution model with six feeder lines has been developed for a local utility network in an urban area. The summer peak load is 23,260 kW which consists of 1,902 customers from categories of residential, small commercial, and industrial (including large commercial). For security reasons, the local utility could not provide detailed information about the network topology and configurations. Therefore, the distribution network is modeled based on the rational alignment of the electrical system and statistical analysis of available data from the local utility and Open Energy Information (OpenEI) dataset [18]. First, the whole system is divided into six sub-regions based on the area map and electrical network topology. Next, we identify a collection of buildings according to the customer demand data collected from a substation of the local utility company. Then the precise number of buildings of each category in each sub-region is derived using the GIS information (Fig. 1). More details on the distribution network modeling can be found in [13].

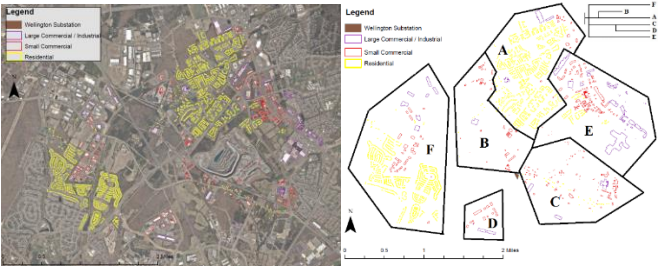

Fig. 1. Study area and sub-regions.

### A. Solar radiation and PV installation

Solar insolation is determined using the LiDAR elevation source data from the US Geological Survey (USGS) "The National Map" (TNM) Download Manager service (U.S. Geological Survey, 2016) which were converted by ArcGIS into a solar insolation raster. Using this process the percentage of each building covered with high insolation points (i.e. average solar insolation of greater than 4.6 kWh/m$^2$/ day) is calculated. Then the buildings will be divided into three categories based on their density of high-insolation coverage. Figure 2 shows the distribution of low, medium, and high-insolation buildings within a selected neighborhood.

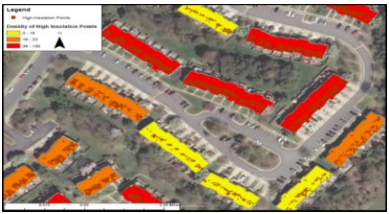

Fig. 2. Distribution of low, medium, and high insolation on buildings.

Next we need to estimate the potential installation capacity of distributed PV in the study area therefore the potential electricity production. First it is assumed that PVDG could only be installed on the rooftops of buildings that fall into the "high-insolation" category. This implies that 510 of the residential buildings (36%), 119 of the small commercial buildings (30%), and 26 of the large commercial and industrial buildings (34%) would become eligible for PV installations. Aggregating the potential PV generation from every eligible building we may determine the total solar power in the studied substation area, It is found that the total potential PV generated power of 16,280 kW is equal to 70% of the total area peak load. The potential annual energy production from PV generation is 21,137 MWh, equal to 18% of the area's annual energy demand of 114,758 MWh. Our proposed study indicates that a substantial portion of the study area's electricity demands could be met through local distributed PV generation.

Considering geographic location of high-insolation buildings we are able to allocate them to a number of buses in the developed distribution network which are potentially ready for a PVDG installation. Consequently, 50 solar ready buses will be considered in the distribution network. Maximum capacity to PVDG installation at each bus is derived based on the number of high insolation buildings connected to the bus and their potential PVDG capacity.

## III. Problem formulation and analysis methodology

To perform a comparative study on the PVDG deployment, we propose two different types of impact assessments. First, an optimization framework is defined to determine the optimal placement and sizing of PVDG systems to manage the power loss and voltage deviation. The objective function is subject to distributed PV constraints and operational constraints of a distribution network, such as avoiding reverse power flows for a given PVDG penetration ratios. Then a stochastic framework is developed to model random PV installations which mimic customer-based renewable integration.

The PV penetration ratio is defined based on the system peak load as follows

$$\gamma(\%) = \frac{\sum_{i=1}^{N} P_{PV_i}^{\max}}{\max_t \left( \sum_{j=1}^{M} P_{Load_j}(t) \right)} * 100\% \quad (1)$$

where $P_{PV}$ and $P_{Load}$ are PV panel output power (kW) and electrical load demand (kW), respectively. The total real energy loss of radial distribution system can be calculated as [19]

$$E_{loss} = \sum_{t=1}^{T} \sum_{l=1}^{L} |i_L^t|^2 R_L \quad (2)$$

where $i_L^t$ is current flowing through line $L$ at time $t$ and $R_L$ is resistance of line $L$. The formulation for voltage profile improvement ($v_D$) with $v_n^t$ as the voltage of bus $N$ at time $t$ is as follow [20]

$$v_D = \frac{1}{TN} \sum_{t=1}^{T} \sum_{l=1}^{N} |v_i^t - 1| \quad (3)$$

### A. Optimal PVDG installation

The multi-objective function can be formulated as follows,

which defines the optimal PVDG siting and sizing for minimizing energy loss and enhancing loadability and voltage profiles while satisfying the network constraints.

$$\min_{\mathcal{L}_{PV}, P_{PV}^{max}} f = E_{loss} + \omega * v_D \quad (4)$$

Subject to:
$$f(P_L, P_{PV}, v|Y_{bus}) = 0 \quad (5)$$
$$i = h(v|Y_{bus}) \quad (6)$$
$$P_{PV}(t) = g(P_{PV}^{max}, I(t)) \quad (7)$$
$$0 \leq P_i^{max} \leq \widehat{P_{PV,i}^{max}} \quad (8)$$
$$1^T P_{PV}^{max} \leq \gamma \cdot P_{sub}^{max} \quad (9)$$
$$i_{ij}^t \geq 0 \quad \forall i < j \quad (10)$$
$$0.95 \leq |v_i| \leq 1.05 \quad (11)$$

where $\mathcal{L}_{PV} = [\ell_1, \ell_2, ..., \ell_n]^T$ $\ell_i \in \{0, 1\}$ is the PVDG location vector, and $\omega$ is a relative weight factor between the two objectives. With $\omega = 0$, the aforementioned optimization problem is equal to minimizing energy loss only. With $\omega = \infty$, it is equal to minimization of voltage deviation alone. $P_{PV}^{max} = [P_1^{PV}, P_2^{PV}, ..., P_n^{PV}]$ is the PVDG installation capacity vector, $Y_{bus}$ is the network admittance matrix, $i$ is vector of bus injected current, $v$ is bus voltage vector, $I(t)$ is solar insolation at time $t$, $\widehat{P_{PV,i}^{max}}$ is PV installation limit for bus $i$ derived from solar data analysis, and $i_{ij}^t$ denotes the current flowing from bus $i$ to $j$ at time $t$, $L$ is the number of lines, $n$ is the total number of buses, and $P_{sub}^{max}$ is defined as the substation peak load. Note that eq. (5) and (6) define the network constraints enforced by AC power flow and network operation constraints, respectively. And eq. (7-9) represent the PVDG installation constraints.

*B. Customer-based integration modeling*

The customer-based integration modeling consists of random siting and sizing of PVDGs which simulate customer decisions on PVDG installation and size selection. With each selected set of location and size of PVDGs the hourly profile of PV generations will be calculated accordingly and fed into the AC power flow model to determine the system state variables, i.e., the bus voltage magnitude and the phase angle. A set of Monte-Carlo experiments will be designed as follows to evaluate the impacts of randomized PV installation on the distribution network operation:

1. Random selection of $S$ locations for PVDG installation from the predefined solar ready buses in the system. This step generates a binary decision vector $= [0/1, ..., 0/1]_{1 \times 50}$, where 1 represents PVDG installation and 0 none-PV installation on the corresponding solar ready bus, with the constraint of $1^T X = S$.
2. Define the PV installation threshold for solar ready buses derived from solar data analysis, i.e., $y = [\widehat{P_{PV}^{max}}]$.
3. Generate the PV size selection factor ($\beta$) using the uniform distribution, $\beta \sim \text{Uniform}[\beta_{\min}, 1]$, where $\beta_{\min} \geq 0$ is called the customer decision factor (CDF) denoting the willingness of customer to install the largest possible PV generation on the site.
4. Determiner the PV installation size for all the selected buses as ( $Z_i = X_i * y_i * \beta_i$).
5. Calculate the PV power output using the solar insolation data and the size of the corresponding PV installation on the site [21].
6. Run a daily time-series AC power flow analysis. The solution results (i.e., voltage, current, and reverse power flow) are stored for the next-step impact assessments.

In this study the local solar insolation profiles have been obtained from [22]. It is worth noting that the customer decision factor (CDF) is set to mimic customer's decision on the PVDG size selection. This factor may be related with various parameters such as the finance budget, incentives, and economics. With CDF we wish to model the willingness or tendency of customers to install high PVDG sizes. A larger value of CDF implies a higher possibility for the customer to utilize all the potential rooftop area to install a largest possible PVDG.

IV. SIMULATION RESULTS AND IMPACT COMPARISON

In this section, the simulation results of optimal vs. random PVDC installations are presented to compare their impacts on the distribution network operations. Particle swarm optimization algorithm is used to solve the optimization problem [23]. In the optimal installation model, we consider the penetration ratios ($\gamma$), ranging from 0% to 70% with 5% step; and three different objective functions: (a) voltage improvement and energy loss reduction; b) Energy loss reduction alone with $\omega = 0$ in (4), and III) Voltage improvement only with $\omega = \infty$.

The algorithm to mimic customer-based random installation is applied to the same distribution network in order to examine the impacts of PVDGs on the distribution network in term of reverse power flow, voltage deviation and energy loss. To generate random PVDG installation samples, customer decision factor is set to ($\beta_{\min} = 0.8$), so that the size selection factor ($\beta$) is a uniform random value between 0.8 and 1.0.

The total reverse power flow experienced by feeders in radial distribution network can be calculated as:

$$F_r^{tot} = \sum_{t=1}^{T} \sum_{l=1}^{L} F_{r,l}^t \quad (12)$$

where $F_{r,l}^t$ denotes the power flow of line $l$ flow at the reverse direction (i.e., feeding back toward the substation) at time $t$. Figure 3 presents a scatter plot of total reverse power flow $F_r^{tot}$ of each random installation case in the given distribution network with the average $F_r^{tot}$ of all the cases with the same penetration ratio depicted as a red solid line. It shows that when the penetration ratio increases, the number of cases that have reverse power flow increases as well. However, the increasing trend does not grow linearly. That is, when the penetration ratio is small i.e. $\gamma < 30\%$, there is not a considerable reverse power flow in random installations. But after some particular penetration ratio, e.g., $\gamma = 30\%$ in our simulated system, a significant raise in the number of cases with reverse power flow issues will appear and more cases will have large magnitudes of $F_r^{tot}$. Note that the distribution network we studies has limited tolerance to operate normally with reverse power flows as reported in [24]. Therefore, the distribution network owners (DNOs) may need to specify a safety threshold for PVDG installation. When the penetration ratio of PVDG installation grows beyond this threshold, mitigation actions such as equipment upgrades or optimized installation of PVDGs will become necessary to manage the reverse power flow issues.

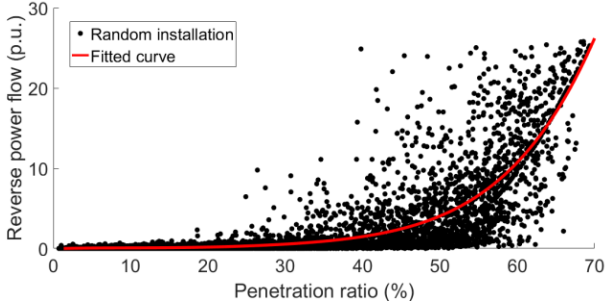

Fig. 3. Total reverse power flow of random installations with $\beta_{min}= 0.8$

The results for voltage deviation across the whole distribution network obtained from both optimized and random installations are presented in Figure 4. For random installations the polynomial function for curve fitting of voltage deviation with respect to penetration ratio is given as:

$$\Delta v = 0.012\gamma^2 - 0.01\gamma + 0.023 \qquad (13)$$

For both optimal installation and the fitting curve of random installations, the voltage deviation achieves its minimum value around the penetration ratio of $\gamma = 40\%$. It also indicates that installing PVDGs with growing penetration ratios will keep improving the system's voltage deviation until it reaches some specific penetration ratio after which the voltage deviation improvement will be decreasing or even diminishing. In comparison with the original system, for the given distribution network, installation of PVDGs improves voltage deviation across the network for all the penetration ratios, which may implies that the DNOs may not need to upgrade their voltage regulatory equipment due to the increasing renewable penetration ratios. However, one must keep in mind that our analysis model in this study has not considered fast dynamics of PV generation caused by cloud moving or stormy weather and the latter may cause severe voltage fluctuations and make necessary the upgrade of voltage regulation scheme.

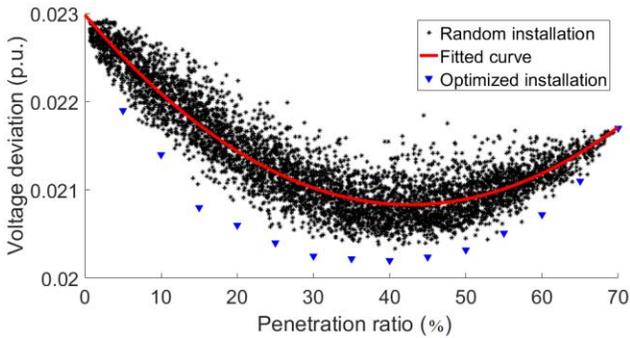

Fig. 4. Voltage deviations of both random and optimized PV installations

Empirical PDFs for voltage deviation of random installations are shown in Figure 5. We observed that the voltage deviations for all penetration ratios follow a normal distribution with different mean values. It can be seen that by increasing penetration ratio mean value for distribution decrease, however, after penetration ratio at 40% mean value increases. In addition, it is found that for high penetration ratios empirical PDF has smaller standard deviation which indicates that for higher penetration ratios the voltage deviation for random samples tend to be close to the mean value, however, for lower penetration ratios (15-25%) the voltage deviations are spread out over a wider range of values.

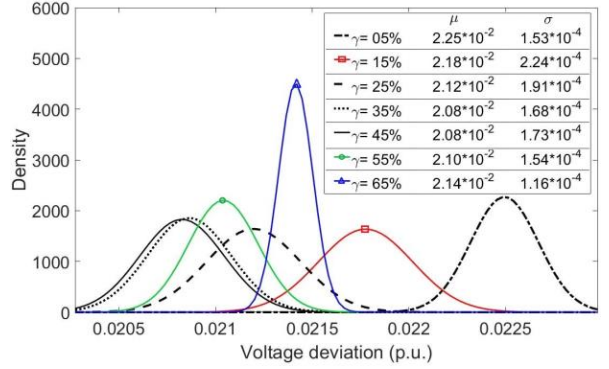

Fig. 5. Empirical PDFs for voltage deviation of random installations

Figure 6 shows the total energy loss of random installations in comparison with optimized installations. For random installations the polynomial function for curve fitting of total energy loss with respect to penetration ratio is as follows:

$$E_{loss} = 3.9 * 10^3 * \gamma^2 - 10^4 * \gamma + 2.5 * 10^4 \qquad (14)$$

It can be seen that by increasing penetration ratio of PVDGs, total energy loss decrease. It is generally accepted that increasing penetration ratio of PVDGs may increase total energy loss in the distribution network for several reasons such as high feeder loadings and lack of local reactive power [25]. However, at given distribution network due to limitation on maximum PVDG penetration forced by available rooftop area for PV installation and solar insolation, penetration ratio does not reach to the critical penetration ratio. Moreover, it is found that there is a significant gap between total energy loss of optimized and random installation of PVDGs in the given distribution network particularly at moderate penetration ratios.

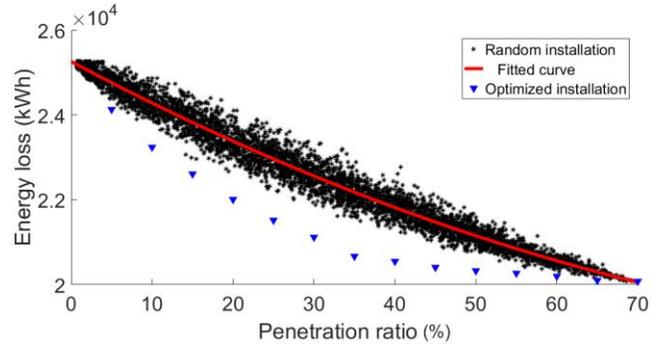

Fig. 6. Energy loss of both random and optimized installations

Empirical PDF for total energy loss of random samples is shown in Figure 7. Our study shows that the total energy loss for all penetration ratios follows a normal distribution with different mean values. It can be observed that as the penetration ratio increase the mean value decrease. Furthermore, increasing the penetration ratio not only changes the mean value but also changes the shape of the empirical PDFs. At penetration ratio of 25% we have a wider PDF (higher standard deviation), however, at higher penetration ratio we have a narrower normal distributions (lower standard deviations). In other words, at higher penetration ratio of random installations, total energy loss of distribution network does not have a large amount of variation.

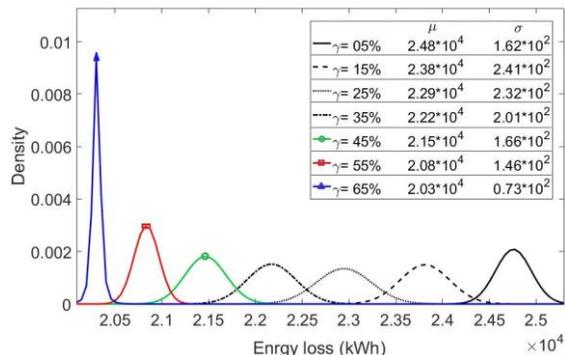

Fig. 7. Empirical PDFs for total energy loss of random installations

## V. Conclusion and Future Works

A comparative study on the PVDG installation in power distribution system is presented in this paper. First, a framework to optimize the siting and sizing of PVDG units is developed with the objective of minimizing the voltage deviation and total energy loss. Then randomized installation of PVDGs is examined to model the customer-based PVDG deployments. Comparing the optimal with the randomized PVDG installations indicates that when the system has a medium renewable penetration ratio an optimal installation is necessary because it will bring significant improvements in energy loss reduction and voltage deviation. However, when the renewable integration ratio is lower or very high, there will be less difference between the types of installation. Depending on DNOs' desire and the expected penetration ratio, the utility company may consider optimized installation at the distribution network. As the future extension of this study we may consider annual analysis along with economic analysis to suggest the best strategy for PVDG installation.